\documentclass[10pt]{article}
\usepackage{amsfonts}
\usepackage{amssymb}
\usepackage{epsfig}
\usepackage{lscape}
\usepackage{comment}
\input{epsf.sty}
\input xypic
\xyoption{all}
\LaTeXdiagrams
\xyoption{2cell}
\UseAllTwocells

\newcommand{\f}[2]{\frac{\displaystyle #1}{\displaystyle #2}}
\def\sq{\sqrt}
\def\sq2{\sqrt{2}}
\def\sq12{\sq{12}}
\def\dsq2{\f{1}{\sqrt{2}}}
\def\be{\begin{equation}}
\def\ee{\end{equation}}
\def\q={\quad = \quad}
\def\lra{\longrightarrow}

\newcounter{examnum}[section]
\newcounter{remarnum}[section]
\begin{document}
\title{The Role of Symmetry in Mathematics}
\author{Noson S. Yanofsky\footnote{Department of Computer and Information Science, Brooklyn College, CUNY. e-mail: noson@sci.brooklyn.cuny.edu}\qquad Mark Zelcer\footnote{Department of Philosophy, SUNY Oswego. e-mail: mark.zelcer@oswego.edu}}
\maketitle
\begin{abstract} 
\noindent Over the past few decades the notion of symmetry has played a major role in physics and in the philosophy of physics. Philosophers have used symmetry to discuss the ontology and seeming objectivity of the laws of physics. We introduce several notions of symmetry in mathematics and explain how they can also be used in resolving different problems in the philosophy of mathematics. We use symmetry to discuss the objectivity of mathematics, the role of mathematical objects, the unreasonable effectiveness of mathematics and the relationship of mathematics to physics. 
\end{abstract}

\section{Introduction}

A\footnote{We wish to thank our friend and mentor Distinguished Professor Rohit Parikh for helpful conversations and for much warm encouragement. Thanks also to Jody Azzouni, Sorin Bangu, Nicolas Fillion, Andr\'{e} Lebel, Jim Lambek, Guisseppe Longo, Jean-Pierre Marquis, Jolly Mathen, Alan Stearns, Andrei Rodin, Mark Steiner, Robert Seely, K. Brad Wray, Gavriel Yarmish, and four anonymous referees, who were all extremely helpful commenting on earlier drafts. N. Y. would also like to thank Jim Cox and Dayton Clark for many stimulating conversations on these topics. He acknowledges support for this project from a PSC-CUNY Award, jointly funded by The Professional Staff Congress and The City University of New York. This work was also supported by a generous ``Physics of Information'' grant from The Foundational Questions Institute (FQXi).} 
 philosophical account of the nature of mathematics must answer at least some of the following questions: First, does mathematics have foundations? If so what are they? If not, why do some branches of mathematics like set theory appear so appealing as foundational accounts? Second, why (and when) are mathematical statements so convincing as compared to other areas of knowledge? Third, how to explain the apparent objectivity of mathematical discourse? Fourth, What is the relationship between mathematics and the sciences? Fifth, Why is mathematics so effective in the world? Sixth, why do the semantics of mathematics appear to match the semantics of natural language? 

The account of mathematics we provide in this paper gives a novel and unified answer to those questions. Our account provides a naturalistic (non-Realist) conception of mathematics that is sensitive to contemporary mathematical practice. We show how we can address the above questions by thinking about mathematics as satisfying a certain set of symmetries. Our account relies on the notions of symmetry that have lately held deep interest for physicists and philosophers of science who begin by inquiring into what it means for laws of physics to satisfy symmetries. Typically, a law satisfies a particular symmetry if it holds despite a change in some parameter of the set of phenomena covered by the law. We will see examples of this and elaborate in Section 2. 

Einstein changed physics forever by taking the idea of symmetry in a novel direction. He showed that rather than looking for symmetries that given laws satisfy, physicists should use symmetries to construct the laws of nature, making symmetry the defining property of laws instead of an accidental feature. These ideas were taken very seriously by particle physicists. Their search for forces and particles are essentially searches for various types of symmetries.  

Philosophers have used these ideas to guide several topics in the philosophy of science. In particular, the symmetries explain the seeming objectivity of the laws of physics. Since the laws of physics are invariant with respect to space, the laws will be the same wherever they are examined. Since the laws of physics are invariant with respect to time, they appear timeless and universal. Some philosophers go further and even question whether the laws of physics actually exist. They posit that rather than there being laws of physics in nature, the physicist is really selecting those families of phenomena that seem to be universal and calling them laws of nature. 

Our goal here is to show that this is true for mathematics as well. The notion of symmetry we focus on we call {\it symmetry of semantics}. Mathematical statements are about abstract or concrete entities of a certain type. Symmetry of semantics says that if the entities are appropriately exchanged for other entities, then the truth value of the mathematical statement remains the same. To a model theorist or a logician this notion is related to validity. Our point is that this type of validity is a form of symmetry. It allows some aspect of the mathematical statement to vary while leaving  its truth value intact. 

Thinking in terms of symmetries will help us deal with different philosophical issues regarding mathematics. We will show that mathematics has no ontological foundations but rather methodological foundations in the creation of symmetries. The apparently ontological foundation of mathematics we see in set theory is merely a special case of a larger kind of mathematical symmetry. The objectivity of mathematical discourse will fall out of the fact that mathematics obeys the symmetries we deemed mathematical. The fact that mathematics works so well in fundamental physics is because the symmetries of physics are a subset of the symmetries of mathematics. And the reason that mathematics appears so convincing is that we determined in advance what counts as mathematics. It is thus not surprising that we are so convinced by the fact that the consequences of our mathematics follow smoothly from our preconceptions of mathematics. 

This paper proceeds as follows: in Section 2 we briefly describe the prominent and evolving role that symmetry plays in the development of physics. In Section 3 we discuss the kinds of symmetries found in mathematics. Section 4 uses these symmetry considerations to address the philosophical questions. Section 5 shows that some of these ideas can be seen within the area of category theory. 

\section{The role of symmetry in physics}

A reasonable history of physics can be given in terms of the ever-expanding place for symmetry in understanding the physical world.\footnote{For philosophical introductions to symmetry see the Introduction in \cite{Brading2003}, \cite{Brading2007}, \cite{Brading2008} and \cite{Bangu2012}. For a popular introduction to the physical issues see \cite{Lederman2004}. For an interesting philosophical view of symmetry in biology see \cite{Bailly2011}. One should also read \cite{LongoForthcoming} for a fascinating perspective about symmetry in modern mathematics.} We briefly outline this expanding history from its origins in the classification of crystals to its role in explaining and defining all the ``laws of nature.'' 

``Symmetry'' was initially employed in science as it is in everyday language. Bilateral symmetry, for example, is the property an object has if it would look the same when the left and the right sides are swapped. In general an object has symmetry if it appears the same when viewed from different perspectives. A cube thus has six-sided symmetry while a sphere is perfectly symmetrical because it looks the same from any of its infinite positions. 

Pierre Curie formulated one of the earliest symmetry rules about nature. He showed that if a cause has a certain symmetry, then the effect will have a corresponding symmetry. Although Curie was mostly concerned with crystals and the forces that create them, the ``Curie Symmetry Principle'' has been employed throughout physics. 

Physicists have generalized the term ``symmetry'' from descriptions of objects to descriptions of laws of nature. A law of nature exhibits symmetry when we can transform the phenomena under the scope of the law in certain ways and still make use of the same law to get a correct result. We say that such a law is ``invariant'' with respect to those transformations. 

A simple example of a symmetry exhibited by a law of nature is the fact that the results of an experiment remain correct when the location of an experiment is changed. A ball can be dropped in Pisa or in Princeton and the time needed for the ball to hit the ground will be the same (all other relevant factors being equal). We thus say that the law of gravity is invariant with respect to location. This fact about locations of experiments was so obvious and taken for granted that scientists did not notice or articulate it as a type of symmetry for some time. 

Similarly, the laws of nature are invariant with respect to time. A physical process can be studied today or tomorrow and the results will be the same. The orientation of an experiment is also irrelevant and we will get identical results (again, {\em ceteris paribus}) regardless of which way the experiment is facing. 

Classical laws of motion formulated by Galileo and Newton display more sophisticated symmetries called ``Galilean invariance'' or ``Galilean relativity.'' These symmetries show that the laws of motion are invariant with respect to an observation in an inertial frame of reference --- the laws of motion remain unchanged if an object is observed while stationary or moving in a uniform, constant velocity. Galileo \cite{Galilei1953} elegantly illustrates these invariances describing experiments that can be performed inside a closed ship. Mathematicians have formulated this invariance by studying different transformations of reference frames. A ``Galilean transformation'' is thus a change from one frame of reference to another that differs by a constant velocity: $x \lra x+vt$, where $v$ is a constant velocity. Mathematicians realized that all these Galilean transformations form a group which they called the ``Galilean group.'' In broad terms then, the laws of classical physics are invariant with respect to the Galilean group. 

To see what it means for laws of physics to be invariant with respect to Galilean transformations consider the following. Imagine a passenger in a car traveling a steady 50 miles per hour along a straight line (neither accelerating nor decelerating). A passenger is throwing a ball up and catching it when it comes down. To the passenger (and any other observer in the moving car) the ball is going straight up and straight down. However to an observer standing still on the sidewalk, the ball leaves the passenger's hand and is caught by the passenger's hand but it does not go straight up and down. Rather the ball travels along a parabola because the ball goes up and down while the car is also moving forward. The two observers are not observing different laws of physics in action. The same law of physics applies to them both. But when the stationary observer sees the ball leaving the passenger's hand, it does not only have a vertical component. Since the car is also moving forward at 50 miles per hour the ball's motion has a horizontal component too. The two observers see the phenomena from different perspectives but the results of the laws must be the same. Each observer must be able to use the same law of physics to calculate where and when the ball will land despite the fact that they make different observations. Thus, the law is invariant with respect to the ability to swap the two perspectives and still get the same answer. The law is symmetric.  One can view this with the aid of the following commutative diagram: 

$$
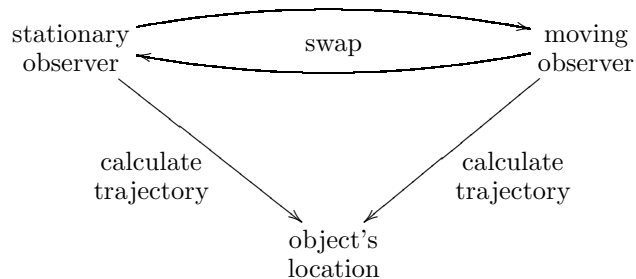
\xymatrix{ 
\txt{stationary \\ observer}  \ar@/^1pc/[rrrr]\ar[ddrr]_{\txt{calculate \\ trajectory}}&&\txt{swap}&& \txt{moving\\observer} \ar@/^1pc/[llll]\ar[ddll]^{\txt{calculate \\ trajectory}}
\\
\\
&& \txt{object's \\location}
}$$
\centerline{\bf Galilean Symmetry}
\vspace{.03 in}
The top part shows the ability of the two perspectives to be swapped. Each observer can calculate the ball's trajectory and both of them must come to the same conclusion about the location of the ball when it lands. 

Another way of expressing this is to say that observers cannot determine whether they are moving at a constant velocity or standing still just by looking at the ball. The laws of physics cannot be used to differentiate between them because the laws operate identically from  either perspective. 

One of the most significant changes in the role of symmetry in physics was Einstein's formulation of the Special Theory of Relativity (STR). When considering Maxwell equations that describe electromagnetic waves Einstein realized that regardless of the velocity of the frame of reference, light will always appear to be traveling at the same rate. Einstein went further with this insight and devised the laws of STR by postulating an invariance: the laws are the same even when the frame of reference is moving close to the speed of light. He found the equations by first assuming the symmetry. Einstein's radical insight was to use symmetry considerations {\em to formulate} laws of physics. 

Einstein's revolutionary step is worth dwelling upon. Before him, physicists took symmetry to be a property of the laws of physics: the laws happened to exhibit symmetries. It was only with Einstein and STR that symmetries were used to characterize relevant physical laws. The symmetries became {\em a priori} constraints on a physical theory. Symmetry in physics thereby went from being an {\em a posteriori} sufficient condition for being a law of nature to an {\em a priori} necessary condition. After Einstein, physicists made observations and picked out those phenomena that remained invariant when the frame of reference was moving close to the speed of light and subsumed them under a law of nature.  In this sense, the physicist acts as a sieve, capturing the invariant phenomena,  describing them under a law of physics, and letting the other phenomena go. 

The General Theory of Relativity (GTR) advanced the relevance of symmetry further by incorporating changes in acceleration. Starting with the advent of GTR Einstein postulated that the laws of nature be understood as invariant even when acceleration is taken into account, i.e. if the observer is accelerating. 

In 1918 symmetry became even more relevant to (the philosophy of) physics when Emmy Noether proved a celebrated theorem that connected symmetry to the conservation laws that permeate physics. The theorem states that for every continuous symmetry of the laws of physics, there must exist a related conservation law. Furthermore, for every conservation law, there must exist a related continuous symmetry. For example, the fact that the laws of physics are invariant with respect to space corresponds to conservation of linear momentum. The law says that within a closed system the total linear momentum will not change. Time invariance corresponds to conservation of energy. Orientation invariance corresponds to conservation of angular momentum, etc (see e.g. \cite{Feynman67} Ch. 4, \cite{Weinberg92} Ch VI, and \cite{Stenger2006} for discussion). Noether's theorem had a profound effect on the workings of physics. Whereas physics formerly first looked for conservation laws, it now looked for different types of symmetries and derived the conservation laws from them. Increasingly, symmetries became the defining factor in physics. 

The ideas of symmetry were significantly advanced by Hermann Weyl. He invented the notion of gauge invariance or gauge symmetry. This is a way of talking about symmetries of laws that are preserved while moving in spacetime. They are a type of local symmetry that moves. This turned out to be very important for almost every branch of physics. 

Researchers currently look for and find some of the more interesting examples of symmetry in particle physics. The field originally postulated three symmetries: a parity invariance with respect to space reflection that lets us swap right and left, a translation from going one way in time to the other, and the charge replacement of a particle with a corresponding anti-particle. Particle physics continues the effort to find more and more abstract symmetries.  The idea is to allow the laws of physics to remain the same no matter how the phenomena are described. 

The physicist Victor Stenger unites the many different types of symmetries under what he calls ``point of view invariance.'' That is, {\em all} the laws of physics must remain the same regardless of how they are viewed. Stenger  (\cite{Stenger2006}) demonstrates how much of modern physics can be recast as laws that satisfy point of view invariance. We can visualize this with a generalization of the previous commutative diagram. 

$$
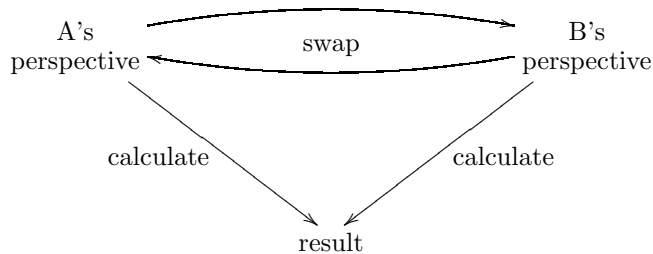
\xymatrix{ 
\txt{A's\\ perspective}  \ar@/^1pc/[rrrr]\ar[ddrr]_{\txt{calculate}}&&\txt{swap}&& \txt{B's \\ perspective} \ar@/^1pc/[llll]\ar[ddll]^{\txt{calculate}}
\\
\\
&& \txt{result}
}$$
\centerline{\bf Point of View Symmetry}
\vspace{.03 in} 
The top part shows the ability of the two perspectives to be swapped. Each perspective can be used to calculate the process of the physical phenomena and both must get the same result. 

Symmetry also plays a role in more speculative areas of physics. 
Our best way forward beyond the standard model are attempts to unify all interactions in nature. One of them, supersymmetry, postulates that there is a symmetry that relates matter to forces in nature. Supersymmetry requires us to postulate the existence of a partner matter particle for every known particle that carries a force, and a force particle for every matter particle. The idea here is that the laws of physics are invariant if we swap all the matter for all the force. None of the partner particles have yet been discovered, but because they are mandated by the symmetries it is what scientists are seeking out. 

Symmetry, as we have described it, is only part of the story. In numerous cases a law of physics actually violates a symmetry law and breaks into several different laws via a mechanism known as ``symmetry breaking.'' These broken symmetries are as conceptually important as the symmetries themselves.  The way a symmetry breaks determines certain constants of nature. But the question of why a symmetry should break in one way and not another is not presently understood. Researchers are at a loss when they leave the constraints set by symmetry.

Recent excitement over the discovery of the Higgs boson reveals a triumph of the role of symmetry in physics. Scientists postulated that there was a symmetry in place at the time of the big bang, and it was only when this symmetry was ``broken'' via the ``Higgs mechanism'' that it was possible for mass to exist. By discovering the Higgs boson, physics was able to provide the mechanism by which mass was produced out of the perfect symmetry of the initial state of the universe. The Higgs mechanism was postulated only on the strength of the presumed symmetries.\footnote{See \cite{Bangu2008} for a related discussion of the discovery of the $\Omega^{-}$ particle.} The recent discovery of the Higgs boson as the culmination of an extensive research program has further vindicated the methodology of postulating symmetries to discover fundamental properties of the universe. 

Physics also respects another symmetry, which as far as we know has not been articulated as such. The symmetry we refer to is similar to the symmetry of time and place that was obvious for millennia but not articulated until the last century. Namely, a law of physics is applicable to a class of physical objects such that one can exchange one physical object of the appropriate type for another of that type with the law remaining the same. Consider classical mechanics. The laws for classical mechanics work for all medium sized objects not moving close to the speed of light. In other words, if a law works for an apple, the law will also work for a moon. Quantities like size and distance must be accounted for, but when a law is stated in its correct form, all the different possibilities for the physical entities are clear, and the law works for all of them. We shall call this invariance for a law of nature {\it symmetry of applicability}, i.e. a law is invariant with respect to exchanging the objects to which the law is applied. We shall see later that this is very similar to a type of symmetry that is central to mathematics. 

To sum up our main point, the change in the role of symmetry has been revolutionary. Physicists have realized that symmetry is the {\em defining property} of laws of physics. In the past, the ``motto'' was that 

\centerline{\underline{A law of physics respects symmetries.}}

\noindent In contrast, the view since Einstein is: 

\centerline{\underline{That which respects symmetries is a law of physics.}}

\noindent In other words, when looking at the physical phenomena, the physicists picks out those that satisfy certain symmetries and declares those classes of phenomena to be operating under a law of physics. 
Stenger summarizes this view as follows ``\ldots the laws of physics are simply restrictions on the ways physicists may draw the models they use to represent the behavior of matter'' (\cite{Stenger2006}: 8). They are restricted because they must respect symmetries. 
From this perspective, a physicist observing phenomena is not passively taking in the laws of physics. Rather the observer plays an active role. She looks at all phenomena and picks out those that satisfy the requisite symmetries. 

This account explains the seeming objectivity of the laws of physics. In order for a set of phenomena to fall under a single law of physics, it must hold in different places, at different times, be the same from different perspectives, etc. If it does not have this ``universality,'' then it cannot be a law of physics. Since, by definition, laws of nature have these invariances,  they appear independent of human perspective. Symmetry thus became fundamental to the philosophical question of the ontology of laws of physics.\footnote{Whether or not there are laws of nature at all or whether they should be eliminated in favor of symmetries is a matter of considerable controversy among philosophers of science. See van Fraassen (\cite{Fraassen89}) and Earman \cite{Earman2004} for stronger and weaker versions of eliminationist views on this issue. Our account is agnostic about this.} 

\section{The role of symmetry in mathematics}

Consider the following three examples: 

(1) Many millennia ago someone noticed that if five oranges are combined with seven oranges there will be twelve oranges in total. It was also noticed that when five apples are combined with seven apples there is a total of twelve apples. That is, if we substitute apples for oranges the rule remains true. In a leap of abstraction, a primitive mathematician formulated a rule that in effect says $5+7=12$. This last short abstract statement holds for any objects that can be exchanged for oranges. The symbols represent any abstract or real entities such as oranges, apples, or manifolds. A similar commutative diagram can be used to illustrate this. 

$$
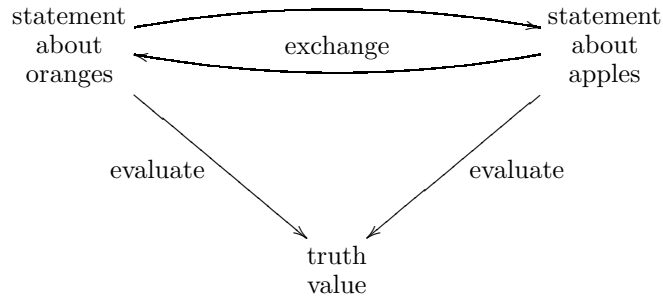
\xymatrix{ 
\txt{statement \\ about \\ oranges }  \ar@/^1pc/[rrrr]\ar[ddrr]_{\txt{evaluate}}&&\txt{exchange}&& \txt{statement \\ about \\ apples} \ar@/^1pc/[llll]\ar[ddll]^{\txt{evaluate}}
\\
\\
&& \txt{truth \\ value}
}$$
\centerline{\bf Symmetry of Fruit Exchange/}
\centerline{\bf Symmetry of Applicability}
\vspace{.03 in} 
The top part shows the ability to swap apples for oranges. Each statement can be evaluated and must produce the same truth value. 

(2) Ancient Egyptians studied different shapes in order to measure the earth so that they can have their taxes and inheritance properly assessed. Their drawings on papyrus represented shapes that could be used to divide up the fields on the banks of the Nile. Archimedes could make the same shapes with sticks in the sand. Drawings can accurately describe properties of these shapes regardless of what they represent: plots of land, Mondrian paintings, shipping containers, or whatever. 

In modern times, mathematicians talk about numerous geometrical or topological theorems such as the Jordan Curve Theorem. This statement says that any non-self-intersecting (simple) closed continuous curve (like a deformed oval) in the plane splits the plane into two regions, an ``inside'' and an ``outside.''  If you exchange one curve for another you will change the two regions. The curve could represent a child's maze or a complicated biological drawing and the Jordan Curve Theorem still applies. 

(3) One of the central theorems in algebra is Hilbert's Nullstellensatz. This says that there is a relationship between ideals in polynomial rings and algebraic sets. The point is that for every ideal, there is a related algebraic set and vice versa. In symbols: $$I(V(J))=\sqrt{J}$$ 
for every ideal $J$. If you swap one ideal for another ideal, you get a different algebraic set. If you change the algebraic set, you get a different ideal. This theorem relates the domains of algebra and geometry and is the foundation of algebraic geometry. 

In these examples we made use of ways of changing the semantics (referent) of mathematical statements. We swapped oranges for apples, changed shapes, transformed curves, and switched ideals. Our central claim is that this ability to alter what a mathematical statement denotes is a fundamental property of mathematics. Of course not all transforms are permitted. If we swap some of the oranges for some of the apples, for example, we will not necessarily get the same true mathematical statement. If we substitute a simple closed curve for a non-simple closed curve (like a figure-8), the Jordan Curve Theorem will not hold true. Such transformations are not legal. We can only change what the statement means in a structured way. Call this structured changing that is permitted a {\em uniform transformation}. 
Our main point is that this uniform transformation and the fact that statements remain true under such a transformation is a type of symmetry. Recall, a symmetry allows us to change or transform an object or ``law'' and still keep some vital property invariant. If a mathematical statement is true, and we uniformly transform the referent of the statement, the statement remains true. Mathematical statements are invariant with respect to uniform transformations. We call the property of mathematical statements that allows it to be invariant under a change of referent {\em symmetry of semantics}. The truth value of the mathematical statement remains the same despite the change of semantic content. 

Symmetry of semantics can be illustrated with the following familiar diagram. 
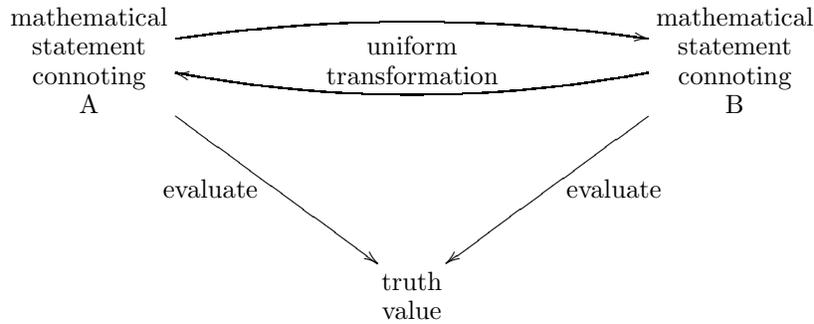
$$\xymatrix{ 
\txt{mathematical\\  statement \\ connoting\\ A }  \ar@/^1pc/[rrrr]\ar[ddrr]_{\txt{evaluate}}&&\txt{uniform \\ transformation}&& \txt{mathematical\\ statement \\connoting\\ B} \ar@/^1.3pc/[llll]\ar[ddll]^{\txt{evaluate}}
\\
\\
&& \txt{truth \\ value}
}$$
\centerline{\bf Symmetry of Semantics}
\vspace{.03 in}
The top part shows the ability to swap connotations of mathematical statements with any two elements of the domain of discourse. Each statement can be evaluated and must arrive at the same truth value. 

Specifically, what types of transformations are uniform transformations? First, the entities we swap must be part of a certain class of elements. Every mathematical statement defines a class of entities which we call its {\em domain of discourse}. This domain contains the entities for which the uniform transformation can occur. When a mathematician says ``For any integer $n \ldots$,'' ``Take a Hausdorff space \ldots'', or ``Let $C$ be a cocommutative coassociative coalgebra with an involution \ldots'' she is defining a domain of discourse. Furthermore, any statement that is true for some element in that domain of discourse is true for any other. A uniform transformation is one in which one element in the domain is substituted for another. Notice that the domain of discourse for a statement can consist of many classes of entities. Each statement might have $n$-tuples of entities, like an algebraically closed field, a polynomial ring and an ideal of that ring. Every mathematical statement has an associated domain of discourse which defines the entities that we can uniformly transform.

Different domains of discourse are indicative of different branches of mathematics. Logic deals with the classes of propositions while topology deals with various subclasses of topological spaces. The theorems of algebraic topology deal with domains of discourses within topological spaces {\em and} algebraic structures. One can (perhaps naively) say that the difference between applied mathematics and pure mathematics is that, in general, the domains of discourse for applied mathematical statements are usually concrete entities while the domains of discourse for pure mathematical statements are generally abstract entities. 

With the concept of domains of discourse in mind one can see how variables are so central to mathematical discourse and why mathematicians from Felix Klein to Tarski, Whitehead (\cite{Epp2011}), Leibniz, Frege, Russell, and Peano all touted their import for mathematics. Variables are placeholders that tell how to uniformly transform referents in statements. Essentially, a variable indicates the type of object that is being operated on within the theory and the way to change its value within the statement. For example in the statement $$a \times (b + c) = (a \times b)+(a\times c)$$ which expresses the fact that multiplication distributes over addition, the $a$ shows up twice on the right side of the equation. If we substitute something for $a$ on the left, then, in order to keep the statement true, that substitution will have to be made twice on the right side. In contrast to $a$, the $b$ and $c$ each occur once on both sides of the equation. Again, the variables show us how to uniformly transform the entities. 

The values of the variables,\footnote{The function of variables (and types) have a long and interesting history. We are not saying that the reason for their creation was to foster the notion of symmetry of semantics. However, variables as they are now, are helpful for dealing with symmetry of semantics. For more on the history of variables and symbols in mathematics see \cite{Mazur2014}, \cite{Heeffer2010}, and \cite{Serfati2005}.} for us, are mathematical objects. They are any entities in a domain of discourse defined by a mathematical statement. So oranges, apples, and stick drawings in the sand are mathematical objects. As long as we can transform those objects into other objects within the same domain of discourse they are mathematical objects. We can transform seven oranges into the elements of the set ${7} = \{0,1,2,3,4,5,6 \}$ and give equal status to each of them as mathematical objects. Mathematicians prefer to use $7$ because of the generality it connotes. But this is misleading. Seven oranges are just as good at representing that number in any mathematical statement. Any statement about the number seven can be made with a transformation of the elements from the set of seven oranges.\footnote{Frege's influence on this definition should be evident. A finite number for Frege consists of the equivalence class of the finite sets where two sets are equivalent if there is an isomorphism from one set to another. When we talk of the equivalence class 5 we are ignorant of the set of the equivalence class under discussion; we may be talking about 5 apples or 5 cars.} The mathematical statements of ``applied mathematics'' are no less true than the statements in ``pure mathematics''. Concrete models of mathematical theories are just as good as abstract models. 

Symmetry of semantics is not a novel concept. It is familiar to logicians and model theorists as the definition of validity. A logical formula is valid if it is true under every interpretation. That is, it must be true for any object in the domain of discourse. The novelty here consists in considering validity as a type of symmetry. We shall see that this symmetry is as fundamental to mathematics as many symmetries are to physics. 

So let us now return to the analogy with physics. Rather than understanding mathematical statements as satisfying symmetry of semantics, we argue that it is that which satisfies these symmetries that we call mathematics. As with physics, in the past we understood that: 

\centerline{\underline {A mathematical statement satisfies symmetry of semantics.}}

\noindent we now claim that:

\centerline{\underline {A statement that satisfies symmetry of semantics is mathematical.}}

\noindent In other words, given the many expressible statements a mathematician finds, her job is to choose and organize those that satisfy symmetry of semantics. In contrast, if a statement is true in one instance but false in another instance, then it is not mathematics. In the same way that the physicist acts as a ``sieve'' and chooses those phenomena that satisfy the required symmetries to codify into physical law, so too the mathematician chooses those statements that satisfy symmetry of semantics and dubs it mathematics. 

Many statements thereby count as mathematical and many do not. Those that do not in general do not satisfy symmetry of semantics; statements containing vague words cannot be mathematical. Even some mathematical-sounding statements such as ``If x is like y and y is like z, then x is like z'' simply fail in most cases because ``like'' is not exact enough to be part of mathematics. 

One may object to this view by saying that allowing mathematics to be whatever satisfies symmetry of semantics is too inclusive. Many general statements not traditionally thought of as mathematical also satisfy symmetry of semantics. For example ``all women are mortal'' is a general non-mathematical statement whereby any woman can be exchanged with any other in the connotation and hence the statement satisfies symmetry of semantics.   
We agree that symmetry of semantics can be found in such general statements, though this is hardly an objection. There is no principled reason to exclude such applications of mathematics.  
Many branches of science, like applied mathematics, make use of general statements with domains of discourse that do not contain ``traditional'' mathematical objects. One would be hard pressed to find a good dividing line between theoretical physics and mathematics, theoretical computer science and mathematics, etc.\footnote{Mark Steiner (\cite{Steiner2005}) treats those fields as applications of mathematics.} Many statements in both pure and applied sciences do satisfy symmetry of semantics and to the extent that a branch of science is mathematical we expect it to have symmetry of semantics. 

We stress again, that many general or universal statements are in fact mathematics. From our perspective if the statement is strong enough so that it is true for {\em every} element of the implied domain of discourse, then it is mathematics. Nor do we shirk away from this definition. Mathematics and its many subdisciplines discuss many different types of objects. No one would say that mathematics is only about numbers and shapes. It is also about propositions, fluid flows, connections on vector bundles, chemical bonds, towels, apples, oranges, etc. Mathematics is about anything where one can reason in an exact manner in such a way that no element in its domain of discourse is exceptional. One might object that by this criterion nothing is outside of our definition of mathematics. This is false. Consider the following statements ``all spoons are silverware'' and ``all silverware is metal'', ``all spoons are metal.'' These are general statements that are not part of mathematics. While it is part of everyday speech --- and would be considered generally true ---  it is not mathematics. Some spoons are not silverware. There are plastic spoons that are neither silverware nor metal. These statements are not exact enough to be part of mathematics. If it was more exact, then the statements would in fact be part of a logic discussion and fall under the dominion of mathematics.\footnote{Exactness comes from the fact that there are no counterexamples. This, in turn, leads to symbolization in mathematics. If we can replace one entity by another, we might as well call the entity $x$.} 

Another symmetry that mathematics has we call {\em symmetry of syntax}. This says that any mathematical object can be described (syntax) in many different ways. For example we can write 6 as $2 \times 3$ or $2+2+2$ or $54/9$. The number $\pi$ can be expressed as $\pi=C/d$, $\pi = 2i\log{1 - i \over 1 +i}$, or the continued fraction $$\pi=3+\textstyle \frac{1}{7+\textstyle \frac{1}{15+\textstyle \frac{1}{1+\textstyle \frac{1}{292+\textstyle \frac{1}{1+\textstyle \frac{1}{1+\textstyle \frac{1}{1+\ddots}}}}}}}$$
Similarly we can talk about a ``non-self-intersecting continuous loop,'' ``a simple closed curve,'' or ``a Jordan curve'' and mean the same thing. The point is that the results of the mathematics will be the same regardless of the syntax we use. Mathematicians often aim to use the simplest syntax possible, so they may write ``6'' or $\pi$ instead of some equivalent statement, but ultimately the choice is one of convenience, as long as each option is expressing the same thing. 

It is interesting to note that formal mathematics started with Euclid's notion of symmetry. The shapes he discussed were essentially invariant with respect to size and orientation. He implored us not to be concerned with how large or small the right triangle was. The important thing was the shape. Similarly for orientation. We are just taking these notions of invariances to be the defining notions of mathematics. We have returned to the origins of formal mathematics. 

The symmetries mentioned above are not the only ones mathematical statements satisfy. Some symmetries are taken for granted to such an extent that even mentioning them seems strange. For example, mathematical truths are invariant with respect to time and space: if they are true now then they will also be true tomorrow, if they are true in Manhattan they are true on Mars. It is similarly irrelevant who asserts a theorem or in what language a theorem is stated, or if it stated at all. 

\section{Some Philosophical Consequences}

{\bf Epistemology.} Considering that we see mathematics as a function of symmetries it is thus natural to think of them as the foundation of mathematics. But let us explore different conceptions of mathematical foundations one at a time. An epistemic foundation for mathematics must explain why mathematical statements are so convincing as compared to other areas of knowledge (\cite{Marquis1995}, \cite{Azzouni2005}). 

On our account, the confidence we have in our mathematics has its origins in our {\em a prioristic} concept of mathematical symmetry. We are certain about mathematical results because we have decided {\em a priori} that the mathematically tractable entities we deal with are those entities that are amenable to what we called uniform transformations. If they are not amenable to such transformations, then we have no reason to be certain they will behave the way we want and so we exclude them from mathematics. The fact that we already decided {\it a priori} how mathematics will work allows us to be certain that our results will turn out the way we expect. We are certain of the outcome because we designed the rules. 

But if certainty in mathematics arises from symmetry considerations and our symmetry considerations are the same as those in science, shouldn't our science provide the same certainty as mathematics? Why are we still less certain about physics than we are about mathematics? The reason is that science and mathematics differ in a crucial way. We generally believe that we have captured all the relevant information in mathematical definitions. In the sciences, it is possible to be unaware of some physical phenomena that ought to be in our domain of discourse or we can fail to understand how some phenomena are captured by a symmetry. This does not happen often in mathematics. To the extent that it does, it mimics the nature of physics and we realize that we should have had a corresponding lack of confidence in our mathematical results. Lakatos' (\cite{Lakatos76}) analysis of the Euler characteristic for the polyhedron is an example of such a case in mathematics. In physics we need look no further than the anomalies that Newton was unaware of or could not handle. They were just not subsumed in his system. Einstein found phenomena that were outside Newton's domain of objects that he could swap in a universal transformation and described an even ``larger'' symmetry that could accommodate them. As long as there are unexplored phenomena (like the perturbations in the perihelion of Mercury or near light speed objects were for Newton) our certainty about science will be lower because we will not know the extent of the domains of discourse that we can uniformly transform or the relevant symmetry of applicability. 

\vspace{2mm} \noindent {\bf Objectivity.} Georg Kreisel is purported to have remarked that the important question in the philosophy of mathematics is the apparent objectivity of mathematical discourse, not the existence of mathematical objects. Mathematical discoveries are sometimes made simultaneously by individuals working independently and the facts of mathematics are true in all places, times, and perspectives. This objectivity is often cast in terms of mathematical realism and has led many to believe in the independent reality of mathematics and its objects. But we can acknowledge the objectivity of mathematics without being realists about mathematical objects by understanding the symmetry at the core of mathematical epistemology. The symmetry we impose reveals that the objectivity of mathematics is an artifact of the way we have designed mathematics. By selecting only those statements that are invariant with regard to what a statement is referring to, the mathematician ensures that the statement is objective and universal. Kant similarly saw mathematical objectivity as a function of the forms of intuitions about space and time; we understand symmetry as the precondition under which mathematics is done. Given a universal precondition, it is unsurprising that we all agree on the truth of the results. We thus need not appeal to an underlying reality for the sake of making sense of objectivity. 

\vspace{2mm} \noindent {\bf Ontology.} The account we have so far given is metaphysically simple. Any object (traditionally mathematical or otherwise) that can be manipulated in a uniform transformation is a mathematical object. We have argued that a mathematical object is any object that is amenable to mathematical treatment. Both seven apples and seven can occupy a domain of discourse. Occupying a domain of discourse is the only relevant criteria for the referent of mathematical statements. That is why it is not odd to assert that ``an apple and an apple are two apples'' and we call that a mathematical statement just as we might call ``$1+1=2$'' a mathematical statement. Moreover, because mathematical discourse appears so much like the discourse of ordinary languages, it has been generally understood since Frege that as a philosophical desideratum we preserve the uniform semantics of our mathematical and ordinary languages (\cite{Benacerraf73}). ``2 is bigger than 1'' has the same structure and truth conditions as ``the Empire State Building is bigger than the Chrysler Building'' and we expect the structural similarities to indicate a similarity in conditions of meaning. And this is exactly what is done when we allow apples, squares, cars, and numbers to be mathematical objects. All can be swapped as part of a uniform transformation so it is unsurprising that our theory of reference looks the same for both natural language and mathematics.

Having a uniform semantics is one thing, but how to handle the ontological question? Prominent accounts of mathematical foundationalism claim that all mathematics can be built out of simpler mathematical stuff. For example some accounts claim that all mathematical objects can be constructed from sets. On the present account, mathematics has no foundation other than the methodological claims of adherence to the symmetry conditions discussed. If that is the case, whence the appeal of such accounts that treat ``foundation'' in connection with mathematics as the ability to show that large parts (or even all) of mathematics can be phrased in some system, and that system is ``primitive.''  Commonly, since many parts of mathematics can be reduced to set theory and logic, and sets provide a convenient domain of discourse in everyday settings, set theory is taken to be a good candidate for such a foundational system. 

This is presumably analogous to the conception of fundamental physics that seeks out the particles in which we, in theory, can express our fundamental ontological statements. This search for fundamental laws or fundamental particles is an important part of contemporary physics. But as we have shown, physics has largely abandoned the idea that programs that search for particles are the starting points for scientific research and theory. Instead we have the presumption that invariances are fundamental, which in turn allow for the discovery of fundamental particles. Thus we really only understand the workings of physics when we look at invariances, not particles. Particles exist, but knowing about them does not give us insight into the nature of the rest of foundational physics. Programs regarding the foundations of mathematics began by confusing reduction and invariance in the same way that physics did before the Einsteinian shift discussed above. By looking at mathematical ``particles'' instead of the invariances needed to understand them we end up in the same situation as pre-Einstein physics; we develop accounts that overlook many things. When those things are discovered the foundational system requires adjustment.  Mathematics is more similar to science than is usually supposed, and thus we must apply the the same rules to both. 

Physics has a number of ``fundamental particles.'' Mathematics admits a variety of types of objects that have no reasonable expectation of reducing one to another (and in many cases it would anyway be unclear what reduces to what). Moreover, we have no one single kind of object with which to phrase all the others. The vocabulary of uniform transformations (i.e. the vocabulary of mathematical methodology) on the other hand is the only way to talk about both abstract and concrete objects. 

Despite the variety of mathematical objects there is still a temptation to claim sets as the foundation. Set theory, after all, allows for the discussion of all kinds of objects. 

However, the fact that we can even reasonably talk about competing foundations for mathematics, exposes a problem. On the present account, the reason sets are not foundational is because of the following: although typical mathematics can be reduced to sets, sets do not exhibit the correct kind of expressive power or display the right kind of symmetries in mathematics to be {\em the} fundamental ``ground''.  To the extent that they appear foundational, sets only display the symmetry of mathematical objects. That is, set theory shows that all mathematical objects are the same in one way: they can all be ``reduced'' to the same thing. Since all mathematical objects are the same there are ways in which they can be treated similarly. This is akin to showing that all (non-fundamental) physical objects reduce to fundamental particles, it fails to deal with of all the other symmetries in nature. As an analogy, recall Frege's definition of a finite number. In it, the equivalence class representing the number 7 happens to contain $\{0, 1, 2, 3, 4, 5, 6\}$. But it also has $\{T, U, P, W, Q, Y, R\}$. Insisting that sets are fundamental to mathematics is akin to insisting that every time we use a set with seven elements  we use $\{0, 1, 2, 3, 4, 5, 6\}$. We can do that but that choice is arbitrary. What is key here is the whole equivalence class and the isomorphisms between the sets. 

Paul Benacerraf (\cite{Benacerraf65}) describes two hypothetical children who are taught about the natural numbers in different ways. Ernie learns that the natural numbers $1, 2, 3 \ldots$ are identified with the sets $\{\varnothing\}, \{\varnothing, \{\varnothing\}\}, \{\varnothing, \{\varnothing, \{\varnothing\}\}\} \ldots$. Johnny learns to identify the natural numbers with the sets $\{\varnothing\}, \{\{\varnothing\}\}, \{\{\{\varnothing\}\}\} \ldots$. The moral of the hypothetical pedagogy is supposed to be that set theory, for example, cannot actually make sense of the myriad of ``fundamental mathematical properties'' because there are an infinite number such of set-theoretical reductions and there is no single set that corresponds to each number.  But on our account this very concern is a symptom of confusing reduction and invariance, not a problem with a particular view of numbers. Benacerraf's problem in other words, is exactly our point. The fact that we can swap $\{\varnothing, \{\varnothing, \{\varnothing\}\}\}$ for $\{\{\{\varnothing\}\}\}$ in a uniform transformation shows that set theory itself exhibits symmetry of semantics. However it says nothing about how sets or set theory are foundational. The fact that we can swap one set for another and understand both as $3$ is not only unsurprising on our account but expected, because set theory is just another branch of mathematics that exhibits symmetry of semantics, like all the others. It is symmetry of semantics that is truly fundamental, not set theory. 

Set theory does not show anything fundamental about numbers because it does not account for how we actually take mathematics to exhibit invariances. Namely, we take mathematical objects to be invariant in a way that the objects stay the same under a wide range of rule transformations, not just object transformations. Therefore set theory initially appears intuitively like a ground for mathematics, but nonetheless fails, because set theory can do one thing that we expect of a physical reduction, namely exhibit something analogous to an ontological-type reduction of some mathematical objects. But the reduction is inadequate as it cannot capture what is really important (what is really mathematical) about mathematics.

Therefore, on this way of looking at mathematics, we need not see any branch of mathematics as ontologically fundamental. The ontology is secondary to what we take to be the methodological underpinnings of mathematics - the search for symmetries.

\vspace{2mm} \noindent {\bf Unreasonable effectiveness.} A philosophical naturalist's interest in the philosophy of mathematics is the alignment of the ontology, epistemology, and especially methodology of mathematics with those of science. The account we have given is naturalistic as it has mathematics relying on the same {\em a priori} role of symmetry as fundamental physics. They both take up the idea that the starting point of inquiry are the postulated symmetries, not the ``smallest pieces.''

Treating mathematics as a function of symmetries also addresses the problem\footnote{Various authors (e.g. Mark Steiner (\cite{Steiner98}) and Nicolas Fillion (\cite{Fillion2012})) now distinguish between various problems of the applicability of mathematics. We confine our remarks to what we take to be Wigner's original question (\cite{Wigner60}) of why mathematics can be used at all with respect to the physical world.} of the unreasonable effectiveness of mathematics in the natural sciences. Eugene Wigner articulated his amazement at the fact that the physical science we discover is shockingly related to the mathematics we need to understand it. Many times science needs to articulate a physical concept and it turns to mathematics; the mathematics is often there. Mark Steiner (\cite{Steiner95}: 154) sees one version of the problem as stemming from the apparent mismatch of methodologies. How can problems emerging from physics be articulated, and even solved, using methods that were designed for a completely unrelated purpose? Another way of looking at it is that if physics and mathematics are both human creations, with physics designed to work with the world and mathematics designed independently of its applications to the natural world, what explains how they work so well together? 

A. Zee, completely independent of our concerns, has re-described the problem as the question of ``the unreasonable effectiveness of symmetry considerations in understanding nature.'' Though our notions of symmetry differ, he comes closest to articulating the way we approach Wigner's problem when he writes that ``Symmetry and mathematics are closely intertwined. Structures heavy with symmetries would also naturally be rich in mathematics'' (\cite{Zee1990}: 319). 

Understanding the role of symmetry however makes the applicability of mathematics to physics not only unsurprising, but expected. Physics discovers some phenomenon and seeks to create a law of nature that subsumes the behavior of that phenomenon. The law must not only encompass the phenomenon but a wide range of phenomena. The range of phenomena that is encompassed defines a set and it is that set which symmetry of applicability operates on. (Recall that symmetry of applicability allows us to exchange one object of a type for another of that type.) So a law must be deliberately designed with symmetry of applicability. Mathematics has a built in ability to express these symmetries because the symmetry of applicability in physics is merely a subset of the symmetry of semantics. That is, the fact that we can exchange one object for another object when dealing with a physical law is simply a special case of exchanging one object for another object in a mathematical statement that expresses the physical law. There is then nothing surprising about the fact that there is some mathematics that is applicable to physics, as the symmetries of physics are a subset of the symmetries of mathematics. Any symmetry we find in physics should (already) be in mathematics. 

For example Newton's established law regarding the relationship of two bodies is $$F=G \frac{m_{1} m_{2}}{r^2} .$$ Symmetry of applicability says that $m_{1}$ can correspond to the mass of an apple or of the moon and the formula still holds. Symmetry of semantics says that $m_{1}$ can be a small number (mass of an apple) or a large number (mass of the moon).  

It is for this reason too that it is odd to say that mathematics is indispensable for physics (in the Quine-Putnam sense). Symmetry of applicability (in physics) is a subset of the mathematical symmetry of semantics. So it is not that mathematics is indispensable for our best scientific theories, but rather, they would not be our best scientific theories (or a recognizable scientific theory at all) if they could not be mathematized.

As we saw, the laws of physics are invariant with respect to the symmetry of applicability. This means that the laws can apply to many different physical entities. Symmetry of applicability is a type of symmetry of semantics. In detail, symmetry of applicability says that a law of nature can apply to many different physical entities of the same type. Symmetry of semantics says that a mathematical statement can refer to many different entities in the same domain of discourse. When a physicist is formulating a law of physics, she will, no doubt, use the language of mathematics to express this law because she wants the law to be as broad as possible. Mathematics shares and increases this broadness. The fact that some of the mathematics could have been formulated long before the law of physics is discovered is not so strange. Both the mathematician and the physicists chose their statements to be applicable in many different contexts. So it is not that mathematics is unreasonably effective, but rather that if it were not effective, it would not be mathematics. The unreasonable usefulness of mathematics is thus no longer mysterious. 

{\bf Mathematical practice}: Our considerations about symmetry would be immaterial if they were not in line with mathematical practice. They are. The day-to-day job of the mathematician is proving theorems. Mathematicians do not generally posit a theorem and then proceed to prove it from axioms. In reality a mathematician has an intuition and formulates some statement. The mathematician tries to prove this statement but almost inevitably finds a counterexample. A counterexample is a breaking (violation) of the symmetry of semantics; there is some element in the supposed domain of discourse for which the statement fails to be true. The mathematician then proceeds to restrict the domain of discourse so that such counterexamples are avoided. Again our indefatigable mathematician tries to prove the theorem but fails, so she weakens the statement. Iterating these procedures over and over eventually leads to a proven theorem. The final theorem may only vaguely resemble the original statement the mathematician wanted to prove. In some sense rather than saying that the ``proof comes to the theorem'' we might say that ``the theorem meets the proof half way.''   
The mathematician acting as a ``sieve'' sorts out those statements that satisfy symmetry of semantics from those that do not, and only those that satisfy this symmetry are reported in the final paper. 
This is just another way of saying that the day-to-day work of the practicing mathematician involves looking for symmetries. Lakatos' ``rational reconstruction'' (\cite{Lakatos76}) may be cited as an example of this constant struggle to preserve the symmetry of semantics of the Euler formula. 

As in physics, we observe symmetry on various levels. We see symmetries being considered when mathematicians construct individual proofs, but also on the level of the formulation of entire mathematical programs. Symmetries in physics are discovered by relying on the intuition that there are substantial domains in which transformations are allowed. It is widely recognized that one way physics progresses is by unification. Unifying an ever larger amount of allowable phenomena under a single given law --- as when Newton united terrestrial and planetary mechanics or Weinberg united the weak force and electromagnetic particle interaction --- is how science advances. Similarly in mathematics; symmetries are discovered when we find that seemingly different mathematical phenomena are really in the same category as an already known transformation and are thereby subsumed under a larger domain; we discover that a new larger class of entities can be uniformly transformed. In other words, we find that there is a union of different domains of discourses which were previously assumed to be comprised of non-interchangeable entities. 

The field of algebraic topology provides an example of such a unification. Researchers realized that there is a certain similarity between taking maps between two topological spaces and taking homomorphisms between two groups. That is, there is a relationship between topological phenomena and algebraic phenomena. Mathematicians went on to use this similarity to try to classify certain topological structures. Category theory grew out of this unification and became a tool for much more unification. Category theory has been derided as ``general abstract nonsense'' that is ``about nothing.'' But precisely because of that, it can be about everything. Hence its language can be used in many different areas of mathematics. The Langlands program is another example of unification. Elegantly described by Frenkel (\cite{Frenkel2013}), the Langlands program is a way of unifying the seemingly different fields of algebraic number theory and  automorphic forms.  As with symmetry, such unification advances mathematics by giving mathematicians an opportunity to discover more general theorems with wider applications and allows them to apply techniques from one domain to the other.\footnote{Philip Kitcher (e.g. \cite{Kitcher76}) touts the importance of these types of cases for mathematics and uses them in the service of demonstrating the existence of mathematical explanation.  Emily Grosholz has studied domain unifications in mathematics extensively. See e.g. \cite{Grosholz00}.} Thus there are various ways in which symmetry considerations aptly describe mathematical practice in the same way they describe scientific practice, lending credence to the idea that this is the proper way to look at mathematics.  

\section{Category theory and the symmetries of mathematics}
An important branch of mathematics that deals with changing objects in mathematical structures is category theory. In this section we discuss the relationship of the notions of symmetry of mathematics with the central notions of category theory. (Knowledge of category theory is not necessary for what follows. Furthermore, this section can be skipped without loss to the philosophical points we made.)
 
A category has objects and morphisms. The objects are usually thought of as mathematical structures and the morphisms are functions between the objects that preserve some aspect of the structure. For example, the category of topological spaces has topological spaces as objects and continuous maps between topological spaces as morphisms. The category of groups has groups as objects and homomorphisms between groups as morphisms. For any algebraic structure, there is a category where the objects are those algebraic structures and the morphisms are functions that preserve some aspect of the structure. In a sense, morphisms are ways of dealing with changing the structure.

One can think of the morphisms in a category as ways of changing some elements in the structure. Let $A$ and $B$ be two objects in a category and let $f: A \longrightarrow B$ be a morphism in the category. In a sense, $a \in A$ gets placed as $f(a) \in B$. The fact that the morphism has to preserve the structure, means that some of the properties of $f(a)$ have to be shared with $a$. Furthermore, some of the properties of $A$ have to be shared with $B$. The strength of the morphism determines what properties of $A$ are in common with properties of $B$. Is the morphism an injection? A surjection? If it is an isomorphism, then $A$ and $B$ have the same categorical properties.
 
One of the central ideas of category theory is that particular constructions are defined by the way morphisms in the category are set up. Most constructions in category theory have ``universal properties'' that describe the construct using morphisms in the category. This is similar to our emphasis that the central idea of a mathematical structure is what is invariant after changing of the elements of the structure. One of the leaders in category theory, F. William Lawvere, summarizes it as follows: ``Thus we seem to have partially demonstrated that even in foundations, not substance but invariant form is the carrier of the relevant mathematical information'' (\cite{LawvereFW:eletcs}).

In this paper we are pushing the notion that one can determine mathematical structures and statements by looking at uniform transformations. There is a very interesting set of ideas in higher algebra and category theory that formalizes this notion of determining structures by looking at uniform transforms. First, some preliminaries. In many places in algebra one looks at an ideal structure and then looks at all the representations / models / algebras of that structure. For example, one can look at 
\begin{itemize}
\item a monoid and the category of sets which the monoid acts on,
\item a group and the category of representations of the group,
\item a ring and its category of modules,
\item a quantum group and its monoidal category of representations,
\item an algebraic theory and its category of algebras,
\item etc.
\end{itemize}
In all these cases one can easily go from the ideal structure to the category of representations. There are times, however, when one can go in the reverse direction. From the category of representations and homomorphisms between representations we can reconstruct the ideal structure. This is similar to the main theme of our paper which is about reconstructing an ideal structure by looking at all the ways objects can be exchanged. 

Exactly how such reconstructions are done is beyond the scope of this paper. However, the core idea of the reconstruction theorems are simple and goes back to the basic definition of what a homomorphism of algebraic structure is. Consider some algebraic structure and let $A$ and $B$ be representations/models/algebras of that algebraic structure. If $+$ is a binary operation, then $f:A \longrightarrow B$ is a homomorphism if the following square commutes:
$$\xymatrix{ 
A\times A \ar[rr]^{+_A} \ar[dd]_{f \times f} && A \ar[dd]^{f} \\
\\
B \times B \ar[rr]_{+_B}&&B.}$$
The reconstructions rest on the idea that the square can be viewed from a slightly different point of view. The usual motto is that 

\centerline{\underline{Homomorphisms are functions that respect all the operations.}}

\noindent We suggest:

\centerline{\underline{Operations are functions that respect all the homomorphisms.}}
\noindent That is, we reconstruct the operations by looking at all those functions that always respect the ways of changing what we are dealing with. If you can swap one element for another and the operation still works, then it is a legitimate operation. In a sense the operations are on equal footing as the uniform transformations.

All these ideas can perhaps be traced back to Felix Klein's Erlangen Program which determines properties of a geometric object by looking at the symmetries of that object. Klein was originally only interested in geometric objects, but mathematicians have taken his ideas in many directions. They look at the automorphism group of many structures, e.g. groups, models of arithmetic, vector spaces, algorithms (\cite{YanofskyG}), etc. In a sense, some of these ideas can be seen as going back to Galois who determined properties of a structure by looking at the set of symmetries of the structure.\footnote{For more about the relationship of category theory and invariance of syntax of semantics see the appendix of \cite{Yanofsky}. See also \cite{Kromer2007} and \cite{marquis2009geometrical} for more on the history and philosophy of category theory and it relationship to symmetry.}  

\bibliography{meirbib,universe}
\bibliographystyle{alpha}
\end{document}